\documentclass[11pt]{article}
\usepackage{amsmath, amsthm, amsbsy, amssymb, graphicx}
\usepackage{adjustbox,url}
\usepackage{subcaption}
\usepackage{pdfsync}
\usepackage{showlabels}
\usepackage[usenames,dvipsnames,table]{xcolor}
\usepackage{algorithm} 
\usepackage{algpseudocode}

\usepackage[hmargin={1in, 1in}, vmargin={1in, 1in}]{geometry}

\usepackage{url, booktabs, wrapfig}
\usepackage{enumitem}
\setlist{noitemsep, topsep=0pt, leftmargin=1em, labelwidth=0em,
  parsep=0pt, partopsep=0pt}

\usepackage{bm,tikz,standalone}

\usepackage{booktabs, threeparttable}

\usepackage{comment} % add comments

\newif\ifshowspecs
\showspecstrue
\showspecsfalse   %%%% uncomment this line in order to remove quoted text from the solicitation
\ifshowspecs
\newcommand{\instructions}[1]{\textcolor{blue}{#1}}
\else
\newcommand{\instructions}[1]{}
\fi

\allowdisplaybreaks
\usepackage{color}
\newcommand{\sid}[1]{{\color{black} #1}}

\usepackage{graphicx}

\usepackage{physics} % for dd

\usepackage[colorlinks=true, citecolor=blue, urlcolor=blue]{hyperref}

\usepackage[]{lineno}
\newcommand*\patchAmsMathEnvironmentForLineno[1]{%
	\expandafter\let\csname old#1\expandafter\endcsname\csname 
	#1\endcsname
	\expandafter\let\csname oldend#1\expandafter\endcsname\csname 
	end#1\endcsname
	\renewenvironment{#1}%
	{\linenomath\csname old#1\endcsname}%
	{\csname oldend#1\endcsname\endlinenomath}}%
\newcommand*\patchBothAmsMathEnvironmentsForLineno[1]{%
	\patchAmsMathEnvironmentForLineno{#1}%
	\patchAmsMathEnvironmentForLineno{#1*}}%
\AtBeginDocument{%
	\patchBothAmsMathEnvironmentsForLineno{equation}%
	\patchBothAmsMathEnvironmentsForLineno{align}%
	\patchBothAmsMathEnvironmentsForLineno{flalign}%
	\patchBothAmsMathEnvironmentsForLineno{alignat}%
	\patchBothAmsMathEnvironmentsForLineno{gather}%
	\patchBothAmsMathEnvironmentsForLineno{multline}%
}

\allowdisplaybreaks % \usepackage{rotating}

\graphicspath{{./}{../figures/}}

\usepackage{lipsum}

\usepackage[compact]{titlesec}
\titleformat{\subparagraph}[runin]
{\normalfont\normalsize\it}{\thesubparagraph}{1em}{}

\usepackage[]{natbib}
\bibpunct{(}{)}{;}{a}{}{,}

\newcommand{\PP}{\mathbb{P}}

\newcommand{\RR}{\mathbb{R}}
\newcommand{\NN}{\mathbb{N}}

\newcommand{\bZ}{\boldsymbol{Z}}

\newcommand{\be}{\boldsymbol{e}}

\newcommand{\bz}{\boldsymbol{z}}

\newcommand{\origin}{\boldsymbol{0}}
\newcommand{\boot}{\text{boot}}

\allowdisplaybreaks

\usetikzlibrary{arrows.meta, positioning, shadows}

\usepackage{authblk}
\title{\textbf{Classification of Extremal Dependence in Financial Markets via Bootstrap Inference}}
\author[1]{Qian Hui}
\author[2]{Sidney I. Resnick}
\author[1,3]{Tiandong Wang \thanks{Corresponding author, \href{mailto:td_wang@fudan.edu.cn}{td\underline{ }wang@fudan.edu.cn}.}}
\affil[1]{Shanghai Center for Mathematical Sciences, Fudan University}
\affil[2]{School of Operations Research and Information Engineering, Cornell University}
\affil[3]{Shanghai Academy of Artificial Intelligence for Science}

\begin{document}
\maketitle

\thanks{{\it In fond memory of Peter Brockwell, a remarkably nice and good
  person who was happy when with his family and friends, happy
  surrounded by
  his books
  and happy in his broadly defined garden.}}
\begin{abstract}
  Accurately identifying the extremal dependence structure in
  multivariate heavy-tailed data is a fundamental yet challenging
  task, particularly in financial applications.  Following a recently
  proposed bootstrap-based testing procedure, we apply the methodology
  to absolute log returns of U.S. S\&P 500 and Chinese A-share
  stocks \sid{over a time period well before the U.S. election in 2024.}
  The procedure reveals \sid{more isolated clustering of dependent assets
  in the U.S. economy compared with China which exhibits different
  characteristics and  a more interconnected pattern of extremal
  dependence.} Cross-market analysis identifies strong
  extremal linkages in sectors such as materials, consumer staples and
  consumer discretionary, highlighting the effectiveness of the
  testing procedure for large-scale empirical applications.
\end{abstract}

\emph{Keywords:} Bootstrap, multivariate regular variation, asymptotic dependence, financial data 

\section{Introduction}
In recent years, financial markets have experienced a number of
extremal events, such as the 2008 global financial crisis and the
COVID-19 pandemic. These events highlight the importance of
understanding the dependence structure between financial assets under
extreme conditions. Traditional correlation-based methods are often
inadequate for capturing the extremal dependence structure, prompting
a shift toward models grounded in extreme value theory. 

For multivariate heavy-tailed data, accurately distinguishing among forms of extremal dependence remains a fundamental and challenging problem. While graphical diagnostics such as angular histograms and scatter plots provide preliminary insights (see for instance \cite{das:resnick:2017}), they often exhibit substantial sensitivity to threshold choice and are insufficient for formal classification. 
Recent advances in the extreme value theory have led to the development of statistical methods for modeling extremal dependence in multivariate settings.
For example, \cite{lehtomaa:resnick:2020} consider the estimation of the support of the angular measure; \cite{hu:peng:sergers:2024} propose a Markov tree-based approach to model multivariate heavy-tailed distributions; \cite{wang:resnick:2025} propose formal hypothesis tests to distinguish different asymptotic dependence structures.
Motivated by the framework introduced in \cite{wang:resnick:2025},
this paper focuses on a systematic and \sid{largely} automated
\sid{classification} procedure for
classifying  asymptotic dependence structures \sid{that allows
categorization of a large number pairs} of financial time
series into four canonical \sid{dependence} types: asymptotic independence, weak
dependence, strong dependence, and full dependence. \sid{These 4
  categories} are characterized by
the support properties of the limit measure arising from \sid{the assumed} multivariate
regular variation. 

As proposed in \cite{wang:resnick:2025}, one way to implement the test
is via bootstrapping, with further justification provided in
\cite{wang:resnick:2025plus}. We apply this testing framework to stock return
data from both the U.S. and China, focusing on the absolute log
returns of selected U.S. S\&P 500 and Chinese A-share stocks. The S\&P
500 is a widely used benchmark index that tracks the performance of
500 large-cap companies listed on U.S. exchanges, representing a broad
cross-section of the U.S. equity market. Similarly, A-shares refer to
stocks of companies based in mainland China, traded on domestic
exchanges such as the Shanghai Stock Exchange (SSE) and the Shenzhen
Stock Exchange (SZSE). Denominated in Chinese yuan, A-shares were
 \sid{restricted to domestic investors until 2003} but have since gained
international visibility through inclusion in global indices like the
MSCI Emerging Markets Index. The selection of S\&P 500 and A-share
stocks offers a comparable basis for examining market dynamics in two
of the world's largest and most influential economies.

For each stock pair, we apply a hierarchical bootstrap-based testing procedure to classify the observed dependence structure. The classification leverages statistical tests sensitive to different configurations of the angular measure: concentration on $\{0,1\}$ (asymptotic independence), on a proper subinterval 
$[a,b]\subsetneq [0,1]$ (strong dependence), on a single point (full dependence), or across the entire interval 
$[0,1]$ (weak dependence). The bootstrap methodology \sid{overcomes
  the difficulty that if parameters are replaced by plug-in
  estimates,then asymptotic distributions become degenerate}
and ensures robustness
to thresholding and estimation variability. 

Empirical results indicate widespread extremal dependence both within
and across the U.S. and Chinese markets. Sectoral differences are
apparent, with U.S. stocks exhibiting \sid{more isolated clustering of
  dependent sectors}, whereas Chinese stocks display more interconnected extremal
behavior. Cross-market analysis reveals sectors such as materials,
consumer staples and consumer discretionary exhibit strong extremal
linkages between the two economies. 

The rest of the paper is organized as follows. Section~\ref{sec:prelim} collects important background knowledge, including the model setup and an overview of
the bootstrap testing framework. Section~\ref{sec:data} gives detailed information about the data analysis, and uncovers different dependence structures in both U.S. and Chinese markets. Concluding remarks are presented in Section~\ref{sec:conclusion}.

\section{Preliminaries}\label{sec:prelim}
\subsection{Dependence Structure}
       Consider iid  data from the common distribution of a random vector
       $\bZ :=(X,Y) \in \RR_+^2$ with $\PP[\bZ \in \cdot ]$
       satisfying
       \begin{equation}\label{e:regVar}
      t \PP[\bZ /b(t)  \in \cdot \,] \to \eta (\cdot),\quad (t\to\infty), 
      \end{equation} 
for a measure $\eta(\cdot)$ on $\RR_+^2\setminus\{\origin\}$ and some regularly varying scaling function $b(t)\to\infty$.
Apply the $L_1$ polar transformation:
$$R=X+Y,\quad \Theta=\frac{X}{X+Y},$$
and the convergence in \eqref{e:regVar} becomes
\begin{equation}\label{e:regVarPolar}
t\PP[(R/b(t) , \Theta) \in \cdot \,] \to \nu_\alpha \times S(\cdot) ,
\quad (t\to \infty),
\end{equation} 
on $(\RR_+\setminus \{0\} )\times [0,1]$ where $\nu_\alpha
(x,\infty)=x^{-\alpha},\,x>0$ and $S(\cdot) $ is a probability measure
on $[0,1]$ 
\sid{called} the \emph{angular measure} \sid{whose support is denoted supp($S$)}. More details can be found in, for example \cite{lindskog:resnick:roy:2014, resnickbook:2024}.

Following \cite{wang:resnick:2025}, we characterize four different cases of asymptotic dependence structures for bivariate heavy-tailed data:
\begin{enumerate}
\item[(i)] Asymptotic Independence: Extremal observations tend to
\sid{lie near} the axes and  $\text{supp}(S)=\{0,1\}$, indicating that
  both components of the data are unlikely to be large
  simultaneously. 
\item[(ii)] Asymptotic Strong Dependence: Extremal observations concentrate within a cone $\mathbb{C}_{a,b}$, where the ratio 
$x/(x+y)$ falls within a specific interval $[a,b]\subsetneq [0,1]$, i.e.  $\text{supp}(S)=[a,b]$.
\item[(iii)] Asymptotic Full Dependence: Extremal observations tend to concentrate on a ray emanating from the origin, i.e. $\text{supp}(S)$ is a single point.
\item[(iv)] Asymptotic Weak Dependence: Extremal observations occur in all directions of the positive quadrant without apparent restrictions, i.e. $\text{supp}(S)=[0,1]$.
\end{enumerate}
For asymptotic independence, we reduce it to
          the full dependence case by the transformation
          $(R,\Theta)\mapsto \bigl(R,g(\Theta)\bigr)$ where $g(0)=g(1)=1$. For
          instance, assume
          $$g(\theta)=\begin{cases}
            1-2\theta,& \text{ if }0\leq \theta < \frac 12,\\
            3-2\theta,& \text{ if }\frac 12 \leq \theta \leq
            1,\end{cases} $$ 
then an asymptotically independent
          heavy tailed distribution is transformed to a distribution with fully
          dependent tail. (Other approaches to handle asymptotic independence
          can be found in
          \citep{lehtomaa:resnick:2020,resnickbook:2024}.)

Assuming an iid sample $\bZ_1,\dots,\bZ_n$ with common distribution
satisfying \eqref{e:regVar}.
For $R_j = X_j+Y_j$, $j=1,\ldots,n$, let $R_{(i)}$ be the $i$th
largest of $R_1,\dots,R_n$, and suppose $\bZ_i^*$ and $\Theta_i^*$ are
the concomitants of $R_{(i)}$. \sid{Extreme value estimation methods
  require choosing either  a threshold or in our case choosing an integer
  $k=k(n)$ representing the number of extreme observations used in estimation.}
Define $d(\bz, \mathbb{C}_{a,b})$ as a
  distance measure of $\bz$ to the cone $\mathbb{C}_{a,b}$:
  for $a>0$,
  \begin{align}
    d((x,y),\mathbb{C}_{a,b})=& \Bigl(\bigl( (b^{-1}-1)x-y\bigr)_+ +
                                \bigl( y-(a^{-1}-1)x\bigr)_+ \Bigr)\label{e:dCart}\\
    \intertext{and in $L_1$-polar coordinates $(r,\theta)=\bigl(x+y,
    x/(x+y)\bigr) $ the distance is}
    =&r\Bigl( \bigl(b^{-1}\theta -1\bigr)_+ +\bigl(1-a^{-1}\theta\bigr)_+\Bigr).\label{e:dPolar}
  \end{align}
If $a=0$, interpret the second term of \eqref{e:dCart} or
  \eqref{e:dPolar} as $0$.

To distinguish the asymptotic dependence structure, \cite{wang:resnick:2025} propose two test statistics
\begin{equation}\label{e:teststats}
          D_n:=\frac{1}{{k(n)}}\sum_{i=1}^{k(n)}
         \left(1+\frac{d\bigl(\bZ_i^*,\mathbb{C}_{a,b}\bigr)}{R_{(k(n))}}\right)\log\frac{R_{(i)}}{R_{(k(n))}},\quad
  T_n :=
          \frac{\sum_{i=1}^{k(n)}\Theta_i^*\log\frac{R_{(i)}}{R_{(k(n))}}}{\sum_{i=1}^{k(n)}\Theta_i^*},
         \end{equation}
and under mild conditions, \cite{wang:resnick:2025} have proved that
two test statistics $D_n=D_n(a,b)$ and $T_n$ are asymptotically normal.
Note that the asymptotic variance of $T_n$ varies by case, and is smallest under full dependence, which helps in classifying these cases.
In addition, for the asymptotic independence case, we apply the $T_n$
statistic to the transformed data $\{\bigl(R_i, g(\Theta_i)\bigr
):1\le i\le n\}$, i.e.
\[
T_n(g) :=\frac{\sum_{i=1}^{k(n)}g(\Theta_i^*)\log\frac{R_{(i)}}{R_{(k(n))}}}{\sum_{i=1}^{k(n)}g(\Theta_i^*)},
\]
which reduces the classification problem to the full dependence case.
Moreover, all statistics require centering by $1/\alpha$ for asymptotic normality, but simply replacing $1/\alpha$ with the traditional Hill estimator 
\citep{hill:1975}, i.e. $\hat{\alpha}=1/\left(\frac{1}{k(n)}
  \sum_{i=1}^{k(n)}\log\frac{R_{(i)}}{R_{(k(n))}}\right)$, causes the
normality results to become degenerate, \sid{necessitating use of the
  bootstrap.}
Theoretical justification of the bootstrap method is provided in \cite{wang:resnick:2025plus}.

\subsection{Methodology}

\begin{figure}[h]
\centering
\begin{tikzpicture}[
    node distance=1.5cm and 2cm,
    base/.style={
        draw, rectangle, rounded corners, 
        font=\footnotesize, inner sep=3pt, 
        drop shadow, align=center
    },
    testbox/.style={
        base, minimum height=1.5cm, minimum width=2.8cm, fill=gray!10
    },
    decisionbox/.style={
        base, minimum height=1cm, minimum width=2.5cm, fill=blue!15
    },
    arrow/.style={->, thick, -{Stealth[length=4pt]}, line width=1pt},
    label/.style={font=\scriptsize, midway, above}
]

% Test nodes
\node (test1) [testbox] {Test $H_0^{(1)}$:\\$S([\hat{a}, \hat{b}]) = 1$\\with $D_n$-statistic};
\node (test2) [below left = 1cm and -.15cm of test1, testbox] {Test $H_0^{(2)}$:\\$\text{supp}(S)$ is a single point \\with $T_n$-statistic};
\node (test3) [below right = 1cm and 0.7cm of test1, testbox] {Test $H_0^{(3)}$:\\$S(\{0,1\}) = 1$\\with $T_n(g)$-statistic};

% Decision nodes
\node (dec_full) [below left= 1cm and -1.5cm of test2, decisionbox] {Asymptotic\\ Full Dependence};
\node (dec_strong) [below right= 1cm and -1.5cm of test2, decisionbox] {Asymptotic\\ Strong Dependence};
\node (dec_asy) [below left = 1cm and -1.5cm of test3, decisionbox] {Asymptotic\\Independence};
\node (dec_weak) [below right = 1cm and -0.5cm of test3, decisionbox] {Asymptotic\\Weak Dependence};

% Arrows
\draw[arrow, thick] (test1.west) -| node[pos=0.8,left] {\scriptsize Accept} (test2);
\draw [arrow, thick] (test1.east) -| node[pos=0.8,right] {\scriptsize Reject}(test3);
\draw [arrow] (test2) -- node[label,left] {Accept} (dec_full);
\draw [arrow] (test2) -- node[label,right] {Reject}(dec_strong);
\draw [arrow] (test3) -- node[label,left] {Accept} (dec_asy);
\draw [arrow] (test3) -- node[label,right] {Reject} (dec_weak);

\end{tikzpicture}
\caption{Overview of the proposed testing procedure.}\label{fig:test}
\end{figure}
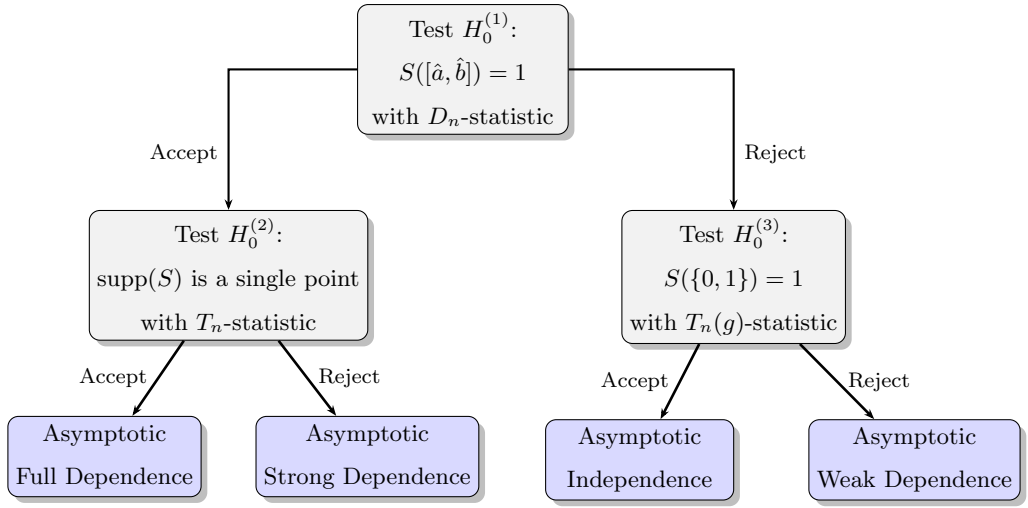

{We \sid{describe} the classification procedure in the following series of
steps which are summarized as a flow chart 
in Figure \ref{fig:test}  and as an algorithm in Algorithm \ref{alg:test}.
Decision steps are as follows:}
\begin{enumerate}
\item[\emph{Step 1.}] For a given sample size $n$, choose
  $k=k(n)\sid{\to \infty}$
  such that $k(n)/n\to 0$, and use the 
            $k(n)$ largest observations for estimates. The bootstrap
            sample size is taken as $m=m(n)\sid{=o(n)}$ following
            \cite{athreya:1987b, gine:zinn:1989, feigin:resnick:1997, {resnickbook:2007}}.
            Also, draw $B$ bootstrap samples of size $m$ from the original
            sample. \sid{These samples are taken independently since they result from multinomial sampling.     }
\item[\emph{Step 2.}] \sid{Using} the original sample, estimate the
              parameter vector $(\alpha, a,b)$. \sid{For $\alpha$ use
                the
    Hill estimator} and $(\hat a,\hat b)$ are \sid{consistently}
    estimated by
    $$(\hat a,\hat b):=\arg\min_{0<a\leq b\leq 1}\Bigl\{(b-a) +
    \lambda \sqrt{k(n)}
    \left\vert D_n - 1/\hat{\alpha}\right\vert \Bigr\},$$
          \sid{\citep[Theorem~5.1]{wang:resnick:2025} where $\lambda>0$
              is a tuning parameter.} When $\lambda$ is large, the
            optimization \sid{penalizes distance to the cone
            and hence tends to } favor a wide estimated interval; we choose $\lambda=4$ in this study.

\item[\emph{Step 3.}]   For a bootstrap
                sample of size $m$, $\{(\Theta_i^\boot, R_i^\boot):1\le i\le m\}$, we write
                \[
                D_m^\boot = \frac{1}{k(m)}\sum_{i=1}^{k(m)}\left(1+\frac{d\bigl(\bZ_i^{*\boot}, \mathbb{C}_{\hat{a},\hat{b}}\bigr)}{R^\boot_{(k(m))}}\right)\log\frac{R^\boot_{(i)}}{R^\boot_{(k(m))}},\qquad 
                T_m^\boot = \frac{\sum_{i=1}^{k(m)}\Theta_i^{*\boot}\log\frac{R^\boot_{(i)}}{R^\boot_{(k(m))}}}{\sum_{i=1}^{k(m)}\Theta_i^{*\boot}}.
                \]
               When there
                are $B$ bootstrap samples we write $D_m^{i,\boot},\,T_m^{i,\boot}$,
                $i=1,\dots, B$.
\item[\emph{Step 4.}] Start by testing the existence of strong dependence: 
$$    H^{(1)}_0:\, S([{\hat a},{\hat b}]) = 1 \text{ vs } H^{(1)}_a:\,
S([{\hat a},{\hat b}])<1,$$
using the simple statistics $D_m^{i,\boot}$,                $i=1,\dots,
B$ and reject if
$$
\left\vert {D}_{m}^{i,\boot}-1/\hat{\alpha}\right\vert > 1.96
\frac{1/\hat{\alpha}}{\sqrt{k(m)}}, \; i=1,\dots,B
$$
is true for at least $5\%$ of the bootstrap samples.
\item[\emph{Step 5.}] Failure to reject $H_0^{(1)}$ 
may be due to $S(\cdot)$ concentrating at \emph{some} single point,
%$\theta \in [\hat a,\hat b]$ 
so we then test for full dependence
$$H^{(2)}_0:\, \text{$\text{supp}(S)$ is a single point}\quad \text{ vs }\quad H^{(2)}_a:
  \text{$\text{supp}(S)$ is not a single point},$$
using $T_m^{i,\boot}, i=1,\dots,m$.
%Since the asymptotic variance of $T_n$ is minimal when full dependence is present, we 
Reject
$H^{(2)}_0$ if the bootstrap samples show excessive variability:
\[
k(m)\frac{ \frac{1}{B-1}\sum_{i=1}^B
   \left(T_{m}^{i,\boot}-\bar{T}_{m}^\boot\right)^2}{1/\hat{\alpha}^2} > \chi^2_{0.95, B-1}/(B-1), 
\]
where $\chi^2_{0.95, B-1}$ denotes the 95\% quantile of a chi-square
distribution with $B-1$ degrees of freedom and $\bar T_m^\boot$ is the
mean of $\{T_m^{i,\boot}, i=1,\dots,m\}$.
\item[\emph{Step 6.}] If $H^{(1)}_0$ is rejected, we further test for asymptotic independence vs weak dependence:
$$H^{(3)}_0: S(\{0,1\}) =1 \quad \text{ vs }\quad
 H^{(3)}_a:  \text{supp}(S) = [0,1],$$
by considering the transformed data $\{(r_i,g(\theta_i)), 1\leq i\leq
n\} $. \sid{If asymptotic independence is present, the transformed
  data possess full dependence, and apply the modified T-statistic}
              \[
              T_{m}^{\boot}(g) := \frac{\sum_{i=1}^{k(m)} g(\Theta^{*\boot}_i) \log\frac{R^\boot_{(i)}}{R^\boot_{(k(m))}}}{\sum_{i=1}^{k(m)} g(\Theta^{*\boot}_i)},
              \]
              Reject $H^{(3)}_0$ if 
\[
k(m)\frac{ \frac{1}{B-1}\sum_{i=1}^B
   \left(T_{m}^{i,\boot}(g)-\bar{T}_{m}^\boot(g)\right)^2}{1/\hat{\alpha}^2} > \chi^2_{0.95, B-1}/(B-1). 
\]           
\end{enumerate}

A schematic summary of the testing procedure is given in Algorithm~\ref{alg:test}.

\begin{algorithm}
  \caption{Testing procedure.}
    \label{alg:test}
	\begin{algorithmic}[1]
	    \Require Estimate $\alpha$, $a$ and $b$ from the original sample, denoted as $\hat{\alpha}$, $\hat{a}$ and $\hat{b}$.
		\State Test $H_0^{(1)}: S([\hat{a},\hat{b}]) = 1$ vs $H_a^{(1)}: S([\hat{a},\hat{b}]) < 1$, using ${D}^\text{boot}_m$.
			\If {Accept $H_0^{(1)}$}
				\State Test $H_0^{(2)}:
                                \text{$\text{supp}(S)$ is a single point}$ vs
                                $H_a^{(2)}: \text{$\text{supp}(S)$ is not a single point}$ using
                                ${T}^\text{boot}_m$. 
					\If {Accept $H_0^{(2)}$}
					\State The random vector $\bZ$ exhibits asymptotic full dependence.
					\Else
					\State The random vector $\bZ$ exhibits asymptotic strong dependence.
					\EndIf
			\ElsIf {Reject $H_0^{(1)}$}
                        \State %{
                        Test $H_0^{(3)}:
                                S(\{0,1\})=1$ vs
                                $H_a^{(3)}:  \text{supp}(S)=[0,1]$,
                                using ${T}^\text{boot}_m(g)$, the
                                modified bootstrap statistic resulting
                                from $\theta\mapsto g(\theta)$.
					\If {Accept $H_0^{(3)}$}
					\State The random vector $\bZ$
                                        exhibits asymptotic
                                        independence. 
					\Else
					\State The random vector $\bZ$ exhibits asymptotic weak dependence.
					\EndIf%}
			\EndIf
	\end{algorithmic} 
\end{algorithm}

\section{Data Analysis}\label{sec:data}
\subsection{Description of Data}
We consider  daily adjusted stock prices of
stocks from the S\&P 500 in the U.S. market and A-shares in the
Chinese market during the period from January 4, 2016, to December 30,
2022. \sid{This time period is prior to the 2024 U.S. election.}
The  66 stocks selected for comparison
represent 11 sectors from both the U.S. and
China, with three stocks per sector.
Table~\ref{tab:stocks} lists the sectors and the representative
  stocks chosen in each sector, and {in Appendix~\ref{sec:append_stocks} we provide a comprehensive 
  review of the chosen stocks across the 11 sectors, including their names, tickers (where applicable), 
  core business operations, and industry classifications.}

Financial \sid{returns} typically exhibit serial dependence, but empirical
studies indicate that lower-frequency sampling (e.g., every other day)
weakens this dependence. \cite{cont:2001} observes that while
financial returns display weak autocorrelation, stronger dependencies
in volatility diminish under coarser sampling intervals. This aligns
with temporal aggregation effects in GARCH models:
\cite{drost:nijman:1993} show that less frequent sampling reduces
volatility persistence. Also, \cite{lo:mackinlay:1988} demonstrate via
variance ratio tests that autocorrelation declines with
lower-frequency data. These studies suggest that coarse sampling
intervals weaken serial dependence. \sid{See also Figure 5 in \cite{wang:resnick:2025}.}

Therefore, to mitigate  serial dependence in stock returns, we compute the
absolute log returns {for each of} these 66 stocks using their
every-other-day 
prices, resulting in a dataset of $n = 822$ observations per stock. We
then apply the bootstrap \sid{classification} method to these absolute returns
{to assess asymptotic dependence between asset pairs. The reasons we chose to
  use
  absolute log returns as opposed to log returns include
  simplicity, reluctance to diminish the data length and because
 prior experience suggests
  that for asset pairs it is rare for a large positive movement to
be matched with a large negative movement.}

\begin{table}[h]
	\centering
	\caption{Sector-based classification of selected stocks from
          U.S. S\&P 500 and Chinese A-Shares in extremal dependence
          analysis.
          %\sid{Suggestion: number the sectors here and then in Figures list
          %numbers; chinese in red, US in blue. This would make the
          %boxes more easily identifiable.}
          }
	\begin{tabular}{llp{7.2cm}}
		\toprule
		Sector& U.S. S\&P 500 & Chinese A-Shares \\
		\midrule
		Communication& GOOG, META, NFLX & Wanda Film, CHINA UNICOM, Phoenix Publishing Media \\
		Energy& CVX, XOM, BP & Shanxi Coking Coal, SINOPEC, SHENERGY \\
		Technology& AAPL, MSFT, AMD & ZTE, IFLYTEK, Inspur Electronic Information \\
		Healthcare& CAH, PFE, BSX & Yunnan Baiyao, Zhifei Biological Products, Aier Eye Hospital \\
		Financials& AMP, BAC, STT & Merchants Securities, Ping An Insurance, Agricultural Bank of China \\
		Consumer Discretionary& AMZN, BKNG, HD & BYD, Midea,
                                                         Hisense Home
                                                         Appliances \\  \\
		Consumer Staples& COST, WMT, TGT & Moutai, Wuliangye, Haitian Flavouring and Food \\
		Industrials& DOV, LMT, EFX & Shenzhen Airport, CRRC, SANY HEAVY INDUS \\
		Materials& LIN, MOS, ECL & BAO IRON, Ganfeng Lithium, WANHUA CHEM \\
		Real Estate& AMT, CSGP, DLR & POLY DEVELOPMENTS, URBAN CONS DEV, Grandjoy \\
		Utilities& AWK, AEE, SRE & WUHAN HOLDING, Dazhong Public, SEP \\
		\bottomrule
	\end{tabular}
	\label{tab:stocks}
\end{table}

{To assess the asymptotic dependence of 66 companies requires a large degree
  of procedural automation in the methodology.}
Initially, we chose the \sid{bootstrap sample} thresholds required by our
  classification procedures by applying the minimum distance
method  \citep{clauset:shalizi:newman:2009, 
virkar:clauset:2014,
  drees:janssen:resnick:wang:2020}
  to obtain $k^*(n)$. 
\sid{Despite knowing the minimum distance method can have
  drawbacks \citep{drees:janssen:resnick:wang:2020},
  the large number datasets meant it was impractical
  to individually analyze each} of them by examining Hill plots to
  check
  how sensible {the minimum distance tail} estimates
  seemed. So we resorted to a formulaic 
  approach that automatically chose the threshold. 
  Mindful of the sample size $n = 822$ for each stock, we
  modified the minimum distance-based threshold as  
\[
80 + \min\{40, (k^*(n)-80)_+\}.
\]

{Using these thresholds, we settled on a Hill estimate for the
  tail index of each stock.   Classifying dependence for each pair of companies
  requires standardization 
  and we applied a power transformation \citep[p. 310]{resnickbook:2007}
    to each pair of stocks so
 that the common tail index was the average of the two individual ones.}
Based on the power transformed absolute log returns, we {then ran 
the classification procedure} between each pair of
stocks. We chose $\lambda = 4$ to estimate $a$ and $b$,
and generated $B=200$ bootstrap resamples with $m=\lceil 6n/k(n)
\rceil$ and $k(m)=\lceil 2m^{0.4} \rceil$. These choices
align with those discussed in \cite{wang:resnick:2025}. 
\sid{Based on experience, we found it useful to make the}
further adjustment that if the estimated interval
length $\hat{b}-\hat{a}\ge 0.85$, then we proceeded to directly test for
asymptotic weak dependence vs asymptotic independence.  
{Otherwise,} we followed the classification procedure outlined in
Algorithm~\ref{alg:test} and Figure~\ref{fig:test}.

Let $\be_i \in \NN^4$ be a four-dimensional vector with the $i$-th
entry equal to 1 and all other entries equal to 0. Each implementation
of the algorithm returns a specific $\be_i$, where $i = 1, 2, 3, 4$
corresponds to one of the four dependence categories: asymptotic
independence, asymptotic weak dependence, asymptotic strong
dependence, and asymptotic full dependence, respectively. The p-value
for each hypothesis test is computed with a significance level set at
$0.025$, where we follow Bonferroni's method to control the type I
error. To ensure stability, the entire bootstrap procedure is repeated
50 times for every pair of stocks, and the average of the resulting
vectors is computed. The averaged vectors are then visualized as shown
in Figures~\ref{fig:stocks-US_China} and \ref{fig:stocks_btw}, where
blue, yellow and gray squares represent asymptotic full, strong and
weak dependence, respectively. \sid{The darker the color}, the more
consistently the 50 repetitions classify the pair into the
corresponding dependence category. 
%\sid{Does this mean
%  white$<$gray$<$yellow$<$blue? Or are there, for example, diffferent
%  intensities of blue that does not show when printed?}

\subsection{Summary of Findings}

Our classification results,
summarized in Figures~\ref{fig:stocks-US_China} and \ref{fig:stocks_btw}, 
{indicate} a lack of asymptotic independence among
U.S. stocks and also among
Chinese stocks. There is also a lack of asymptotic independence
between Chinese and U.S. companies.
This {presumably results from} the fact that the global economic world is
flat \citep{friedman:2005} with 
connected  global 
markets and therefore  extreme events in  one region propagate across
borders; this {presumption} is especially {plausible}
between the two  largest economies in the world.

\subsubsection{Comparing sectoral dependence within the U.S. and Chinese markets}

\begin{figure}
\centering
\includegraphics[width=\textwidth]{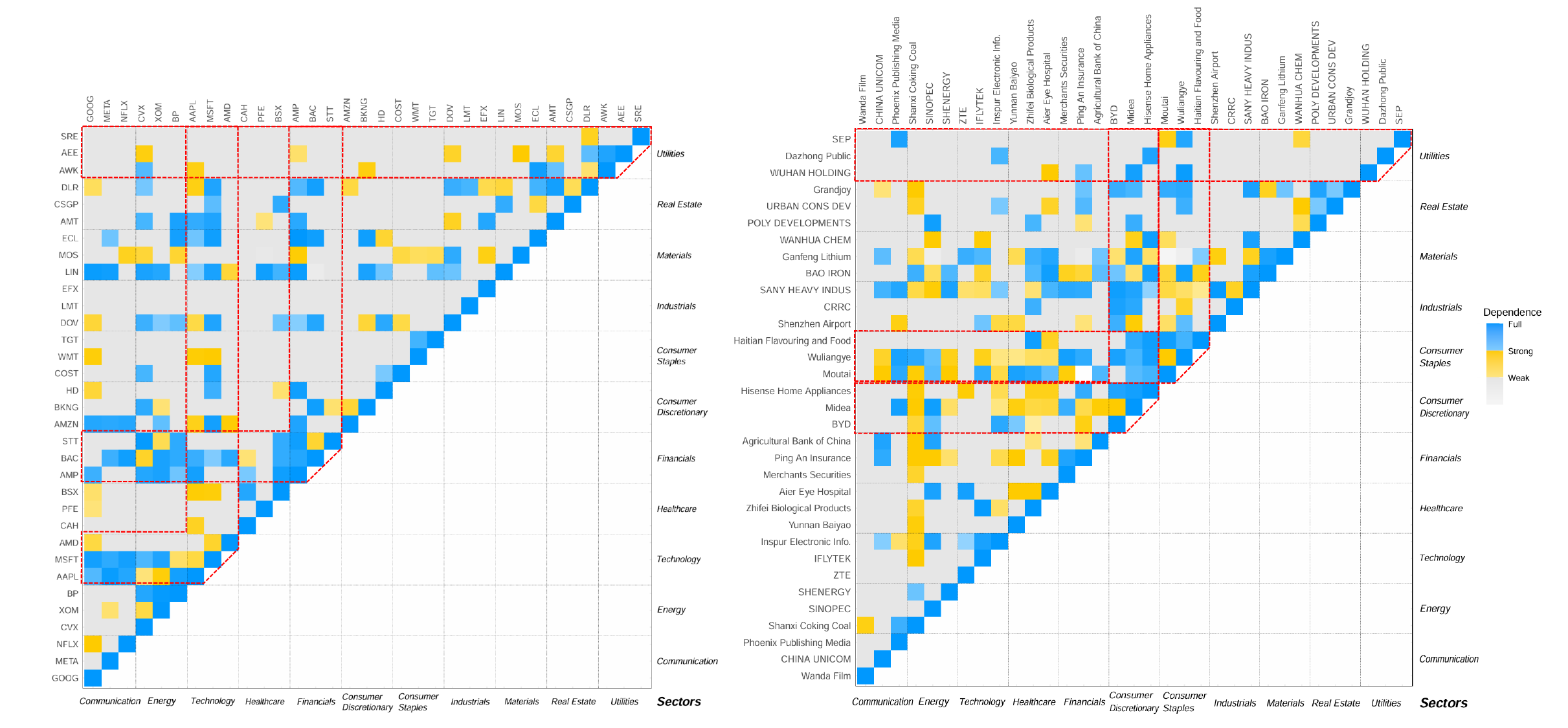}
\caption{Extremal dependence structures; blue, yellow and gray squares
          represent asymptotic full, strong and weak dependence,
          respectively. Left: U.S. market returns segmented by sector.
          Right: Chinese market returns by sector.}
\label{fig:stocks-US_China}
\end{figure}

%\begin{figure}[h]
%\centering
%\includegraphics[width=\textwidth]{Figure/China-f.pdf}
%\caption{Left: Extremal dependence structures within the Chinese market; blue, yellow and gray squares
 %         represent asymptotic full, strong and weak dependence,
  %        respectively.
   %       Right: Extremal dependence structures of Chinese utilities, consumer staples and consumer discretionary sectors as illustrative examples.}\label{fig:stocks-China}
%\end{figure}

\begin{figure}[h]
\centering
\includegraphics[width=\textwidth]{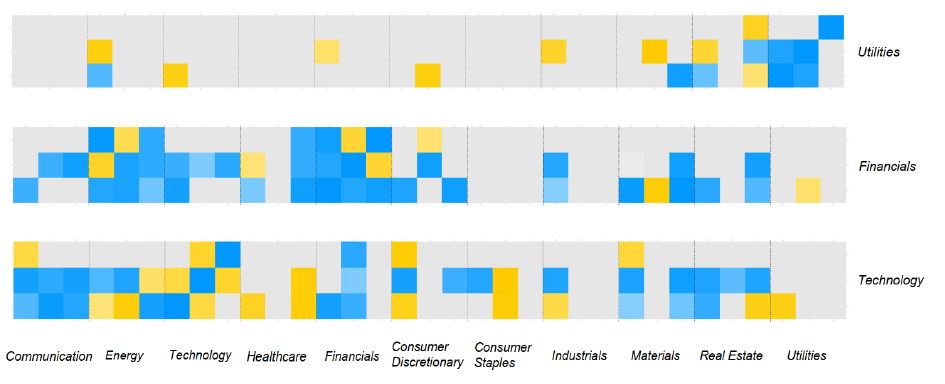}
\caption{Extremal dependence of sectors within the U.S.
  market emphasizing utilities, financials and technology.}\label{fig:US_sector}
\end{figure}

\begin{figure}[h]
\centering
\includegraphics[width=\textwidth]{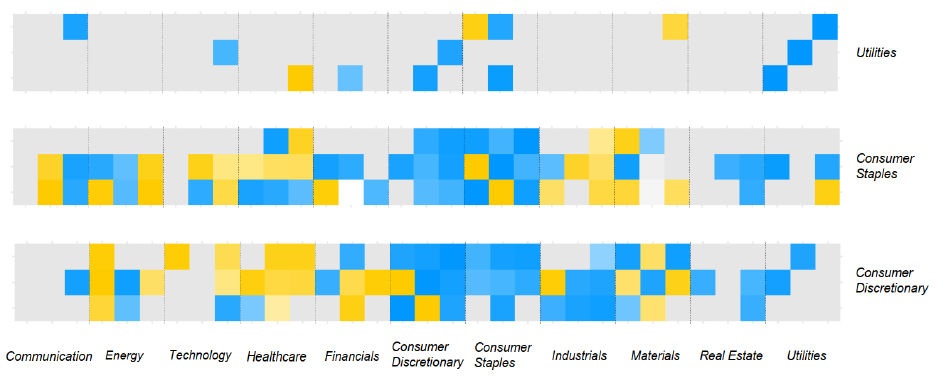}
\caption{Extremal dependence of sectors within the Chinese
  market emphasizing utilities, consumer staples and consumer
  discretionary.}\label{fig:China_sector}
\end{figure}

The left and right panels of Figure~\ref{fig:stocks-US_China}
separately summarize the extremal dependence structure between
\sid{the chosen stocks for the}
economic sectors for the U.S. market (left) and China (right). Due to
symmetry, we only show the upper triangular portion of the dependence
graphs. 
In the left panel of Figure~\ref{fig:stocks-US_China}, we highlight
(in red) three sectors in the U.S. market: utilities, financials, and
technology. For the utility sector, although the lower triangle of the
last \sid{$3\times 3$}
box is blank, it can be completed by mirroring the upper triangle
along the diagonal, resulting in the first row of
Figure~\ref{fig:US_sector}. Similarly, for the financials and
technology sectors, we mirror the highlighted vertical boxes along the
diagonal to obtain the second and third rows of
Figure~\ref{fig:US_sector}, respectively.  For example, the 6 boxes on
the right of the  ``financials" bar are 
obtained from the U.S. portion of Figure \ref{fig:stocks-US_China} by reading
 down the financials column and then rotating each $3\times 3$ box to account
 for $(j,i)$ vs $(i,j)$.
Following a similar strategy, we highlight the utilities, consumer staples, and consumer discretionary sectors in the Chinese market, with properly
mirrored dependence structure displayed in Figure~\ref{fig:China_sector}.

{Visual comparison of} the two \sid{panels of Figure \ref{fig:stocks-US_China}}
indicates structural
differences in the two countries. 
Compared to the left panel, the right one
visually shows less gray
  area, suggesting the Chinese market has strong\sid{er}
extremal dependence among sectors compared to the
U.S. market (especially for consumer staples and consumer
discretionary) while dependence of sectors in
the U.S. market appear more structured and 
segmented compared with China.
The sectors of the Chinese
economy appear more 
  interdependent. 
  The Chinese economy
  is characterized by a 
rapid pace of development, significant government influence,emphasis
on exports over domestic consumption and more
integrated supply chains. {This results in} a greater degree of extremal
dependence across sectors. Such interconnectedness reflects China's
economic growth model, where sectors often rely on each other for
expansion, and the government plays a key role in managing sectoral
relationships and economic growth. 

Some comments on internal dependencies within each country:
We isolate certain sectors for closer inspection. For the U.S. see 
Figure \ref{fig:US_sector} which is obtained from the left panel
  of Figure \ref{fig:stocks-US_China}.
  Certain sectors of the U.S. economy such as
  energy (not illustrated in Figure \ref{fig:US_sector} but can be
  discerned in Figure~\ref{fig:stocks-US_China}),
technology (bottom row of Figure \ref{fig:US_sector} and financials
(middle row in Figure \ref{fig:US_sector})
{exhibit darker colors} and thus
more dependency. 
Other sectors, e.g. healthcare,  consumer staples {and utilities (see
top bar of Figure \ref{fig:US_sector})}  {show more
gray} and thus are {less dependent} on other parts of the U.S. economy.
This suggests that sectors like financials and technology may play a more central role in driving broader economic dynamics.

For China (see the bottom two bars in Figure \ref{fig:China_sector})
we sees heavy doses of blue for the sectors {\it consumer
  staples\/} and {\it consumer discretionary\/} indicating these
{sectors are  heavily dependent on other economic
  sectors.} Presumably {this reflects Chinese manufacturing expertise
of consumer goods.}

Despite these differences between the U.S. and China, both markets
show some similarities in the utilities sector, which remains less
dependent on other sectors for both countries, though with
differences.  Compared to China, the U.S. utility sector, as revealed
by the top rows of Figures~\ref{fig:US_sector} and \ref{fig:China_sector}, shows
stronger extremal dependence on real estate, followed by slightly
weaker dependence on materials and energy, and
weak dependence on the other sections.  These are possibly due to the close connection
between utilities and property development, as well as the need for
raw materials and energy in infrastructure projects. However,
from the top rows of
Figures~\ref{fig:US_sector} and \ref{fig:China_sector},  we see that compared to the U.S.,
the Chinese utilities sector
exhibits somewhat stronger extremal dependence on the consumer
discretionary and staples sectors, indicating that consumer spending
may play an important role in shaping the utility demand.

\sid{Although the Chinese economy has
traditionally been export-driven, 
our findings are consistent with current Chinese government policy
designed to increase consumer spending.}
The utilities sector's dependence on consumer spending may
reflect an emerging shift towards a more consumption-driven economy as
well as government efforts to stabilize and support domestic
consumption.

\subsubsection{Sectoral dependence between the two economies}

\begin{figure}[h]
    \centering
    \includegraphics[width = \textwidth]{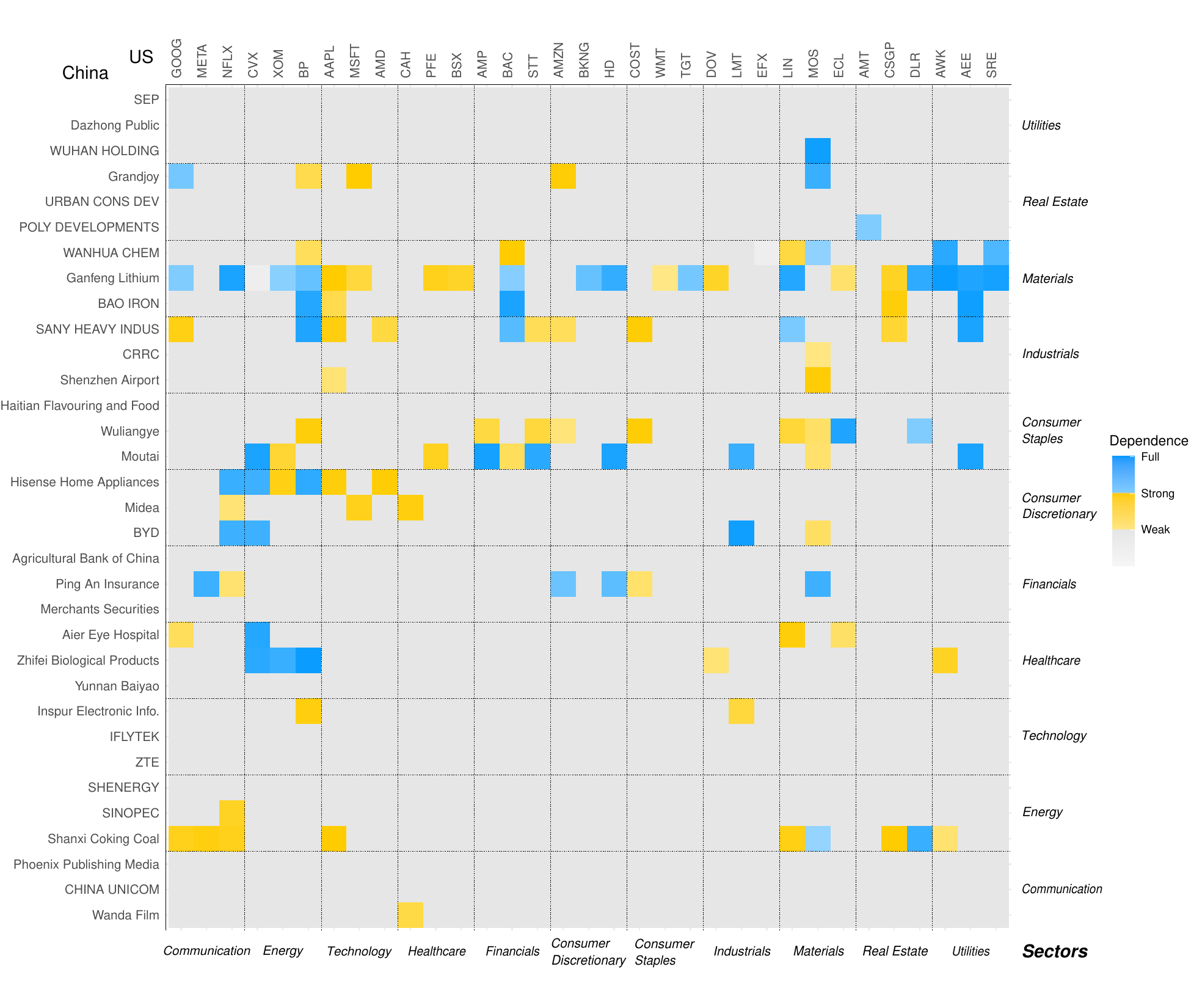}
    \caption{Extremal dependence structures between the U.S. and Chinese stock markets. Blue, yellow and gray squares represent asymptotic full, strong and weak dependence, respectively.}
    \label{fig:stocks_btw}
\end{figure}

Figure~\ref{fig:stocks_btw} displays the extremal dependence
between Chinese and the U.S. stocks grouped by sector, where we find
various patterns of interconnectedness across sectors.  
Between China and the U.S., sectors such as materials, consumer
discretionary, and consumer staples show extremal 
dependence, indicating that extreme market changes in one country
are likely to be reflected in the other. This may be due to deep trade
linkages, global supply chains, and investment flows.

Meanwhile, financials and industrials sectors exhibit weak extremal
dependence {between the two countries}. Although financial markets are
globally linked, China's 
capital controls and state-regulated banking limit direct
spillovers. However, interest rate changes, trade tensions, and
economic downturns may still affect financial institutions in both
countries. The {horizontal bar corresponding to} the
Chinese industrial sector {is also colored predominantly 
  gray} but remains
partially exposed to manufacturing and trading due to government
policies, long-term contracts, and diversified trade partnerships. 

{As Figure \ref{fig:stocks_btw} shows,}
the communication, energy, technology, real estate, and
utilities sectors show little extremal dependence, suggesting a {minimal}
dependence under extreme market conditions. {In China,} communication and
technology are heavily regulated and {in both countries}
geopolitically sensitive, while
energy markets are influenced more by global commodity prices and geopolitical dynamics than by bilateral relations. 
Real estate is shaped largely by domestic policies and housing demand, 
while utilities, as essential services, tend to remain stable and resilient to financial volatility. 

Overall, these findings highlight varying degrees of interdependence, with consumer-driven and materials sectors being the most interconnected, while more regulated or localized industries remain relatively unaffected by economic shocks.

\section{Concluding Remarks}\label{sec:conclusion}
This paper implements a systematic framework for classifying extremal dependence structures in financial time series.
Applied to large-scale stock return data from the U.S. and Chinese
markets, our empirical analysis uncovers both within- and cross-market
asymptotic dependence structures. In particular, the U.S. market
\sid{exhibits stronger clustering of extremally dependent assets}
while
Chinese stocks show more interconnected extremal behavior. Across the
two markets, consumer-related sectors appear to have the strongest
dependence, reflecting underlying economic linkages such as global
supply chains and consumption-driven risk propagation. 

Future work may extend this approach to dynamic settings, incorporate network-based interpretations of dependence, and apply similar classification methods to multivariate extremes in other areas such as climate science, epidemiology and infrastructure risk. The intersection of statistical theory and practical financial modeling continues to offer a rich domain for methodological advances and impactful applications.

\section*{Acknowledgments}
T. Wang gratefully acknowledges Science and Technology Commission
of Shanghai Municipality Grant 23JC1400700 and National Natural Science
Foundation of China Grant 12301660.

T. Wang also thank Shanghai Institute for Mathematics and Interdisciplinary Sciences (SIMIS) for their financial support. This research was
partly funded by SIMIS under grant number SIMIS-ID-2024-WE. T. Wang
is grateful for the resources and facilities provided by SIMIS, which were
essential for the completion of this work.

\bibliographystyle{chicago}      % Chicago style, author-year citations
%\bibliography{/Users/tiandongwang/Documents/bootstrap_fulldep/JTSA/Bootstrap_method/bibfile.bib}
\bibliography{bibfile.bib}

\def\cprime{$'$}
\begin{thebibliography}{}

\bibitem[\protect\citeauthoryear{Athreya}{Athreya}{1987}]{athreya:1987b}
Athreya, K.~B. (1987).
\newblock Bootstrap of the mean in the infinite variance case.
\newblock {\em Ann. Statist.\/}~{\em 15\/}(2), 724--731.

\bibitem[\protect\citeauthoryear{Clauset, Shalizi, and Newman}{Clauset
  et~al.}{2009}]{clauset:shalizi:newman:2009}
Clauset, A., C.~Shalizi, and M.~Newman (2009).
\newblock Power-law distributions in empirical data.
\newblock {\em SIAM Rev.\/}~{\em 51\/}(4), 661--703.

\bibitem[\protect\citeauthoryear{Cont}{Cont}{2001}]{cont:2001}
Cont, R. (2001).
\newblock Empirical properties of asset returns: {S}tylized facts and
  statistical issues.
\newblock {\em Quantitative finance\/}~{\em 1\/}(2), 223.

\bibitem[\protect\citeauthoryear{{Das} and {Resnick}}{{Das} and
  {Resnick}}{2017}]{das:resnick:2017}
{Das}, B. and S.~{Resnick} (2017).
\newblock Hidden regular variation under full and strong asymptotic dependence.
\newblock {\em Extremes\/}~{\em 20\/}(4), 873--904.

\bibitem[\protect\citeauthoryear{Drees, Jan{\ss}en, and {Resnick, S.I., Wang,
  T.}}{Drees et~al.}{2020}]{drees:janssen:resnick:wang:2020}
Drees, H., A.~Jan{\ss}en, and {Resnick, S.I., Wang, T.} (2020).
\newblock On a minimum distance procedure for threshold selection in tail
  analysis.
\newblock {\em Siam J. Math. Data Sci.\/}~{\em 2\/}(1), 75--102.

\bibitem[\protect\citeauthoryear{Drost and Nijman}{Drost and
  Nijman}{1993}]{drost:nijman:1993}
Drost, F. and T.~Nijman (1993).
\newblock Temporal aggregation of {GARCH} processes.
\newblock {\em Econometrica\/}~{\em 61\/}(4), 909--927.

\bibitem[\protect\citeauthoryear{Feigin and Resnick}{Feigin and
  Resnick}{1997}]{feigin:resnick:1997}
Feigin, P. and S.~Resnick (1997).
\newblock Linear programming estimators and bootstrapping for heavy tailed
  phenomena.
\newblock {\em Adv. in Appl. Probab.\/}~{\em 29}, 759--805.

\bibitem[\protect\citeauthoryear{Friedman}{Friedman}{2005}]{friedman:2005}
Friedman, T.~L. (2005).
\newblock {\em The World is Flat: A Brief History of the Twenty-First Century}.
\newblock Macmillan.

\bibitem[\protect\citeauthoryear{Gin{\'e} and Zinn}{Gin{\'e} and
  Zinn}{1989}]{gine:zinn:1989}
Gin{\'e}, E. and J.~Zinn (1989).
\newblock Necessary conditions for the bootstrap of the mean.
\newblock {\em Ann. Statist.\/}~{\em 17\/}(2), 684--691.

\bibitem[\protect\citeauthoryear{Hill}{Hill}{1975}]{hill:1975}
Hill, B. (1975).
\newblock A simple general approach to inference about the tail of a
  distribution.
\newblock {\em Ann. Statist.\/}~{\em 3}, 1163--1174.

\bibitem[\protect\citeauthoryear{Hu, Peng, and Segers}{Hu
  et~al.}{2024}]{hu:peng:sergers:2024}
Hu, S., Z.~Peng, and J.~Segers (2024).
\newblock Modeling multivariate extreme value distributions via {M}arkov trees.
\newblock {\em Scandinavian Journal of Statistics\/}~{\em 51\/}(2), 760--800.

\bibitem[\protect\citeauthoryear{Lehtomaa and Resnick}{Lehtomaa and
  Resnick}{2020}]{lehtomaa:resnick:2020}
Lehtomaa, J. and S.~Resnick (2020).
\newblock Asymptotic independence and support detection techniques for
  heavy-tailed multivariate data.
\newblock {\em Insurance: Mathematics and Economics\/}~{\em 93}, 262 -- 277.

\bibitem[\protect\citeauthoryear{Lindskog, Resnick, and Roy}{Lindskog
  et~al.}{2014}]{lindskog:resnick:roy:2014}
Lindskog, F., S.~Resnick, and J.~Roy (2014).
\newblock Regularly varying measures on metric spaces: Hidden regular variation
  and hidden jumps.
\newblock {\em Probab. Surv.\/}~{\em 11}, 270--314.

\bibitem[\protect\citeauthoryear{Lo and MacKinlay}{Lo and
  MacKinlay}{1988}]{lo:mackinlay:1988}
Lo, A.~W. and A.~C. MacKinlay (1988).
\newblock Stock market prices do not follow random walks: {E}vidence from a
  simple specification test.
\newblock {\em The Review of Financial Studies\/}~{\em 1\/}(1), 41--66.

\bibitem[\protect\citeauthoryear{Resnick}{Resnick}{2007}]{resnickbook:2007}
Resnick, S. (2007).
\newblock {\em Heavy Tail Phenomena: Probabilistic and Statistical Modeling}.
\newblock Springer Series in Operations Research and Financial Engineering. New
  York: Springer-Verlag.
\newblock ISBN: 0-387-24272-4.

\bibitem[\protect\citeauthoryear{Resnick}{Resnick}{2024}]{resnickbook:2024}
Resnick, S. (2024).
\newblock {\em The Art of Finding Hidden Risks; Hidden Regular Variation in the
  21st Century}.
\newblock Switzerland: Springer.
\newblock isbn: 978-3-031-57598-3.

\bibitem[\protect\citeauthoryear{Virkar and Clauset}{Virkar and
  Clauset}{2014}]{virkar:clauset:2014}
Virkar, Y. and A.~Clauset (2014).
\newblock Power-law distributions in binned empirical data.
\newblock {\em Ann. Appl. Stat.\/}~{\em 8\/}(1), 89--119.

\bibitem[\protect\citeauthoryear{Wang and Resnick}{Wang and
  Resnick}{2025a}]{wang:resnick:2025}
Wang, T. and S.~Resnick (2025a).
\newblock Distinguishing forms of asymptotic dependence in heavy tailed data.
\newblock {\em Statistica Sinica\/}.
\newblock To appear; DOI=10.5705/ss.202024.0196.

\bibitem[\protect\citeauthoryear{Wang and Resnick}{Wang and
  Resnick}{2025b}]{wang:resnick:2025plus}
Wang, T. and S.~I. Resnick (2025b).
\newblock Bootstrap methods for testing asymptotic dependence in multivariate
  heavy-tailed data.
\newblock Manuscript in preparation.

\end{thebibliography}

\appendix
\section{Summary of Chosen Stocks}\label{sec:append_stocks}
This summary categorizes companies by industry, listing their names, tickers (where applicable), and core business operations, and all
information is sourced from \href{https://finance.yahoo.com}{Yahoo Finance}.
\subsection*{Communication}
- \textbf{GOOG – Alphabet Inc.}: Offers various products and platforms, including Google Search, YouTube, Android, and Google Cloud, across multiple regions.  \\
- \textbf{META – Meta Platforms, Inc.}: Develops products enabling
people to connect and share through mobile devices, personal
computers, and virtual reality headsets. \sid{Includes Facebook,
  WhatsApp, Instagram.} \\
- \textbf{NFLX – Netflix, Inc.}: Provides entertainment services, offering television series, documentaries, feature films, and games across various genres and languages. \\ 
- \textbf{Wanda Film – Wanda Film Holding Co., Ltd.}: Engages in the investment, construction, and operation of movie theaters in China, Australia, and New Zealand.  \\
- \textbf{CHINA UNICOM – China United Network Communications Limited}: Provides various telecommunication services, including mobile network, broadband, and mobile data services in China.  \\
- \textbf{Phoenix Publishing Media – Jiangsu Phoenix Publishing \& Media Corporation Limited}: Engages in editing, publishing, and distribution of books, newspapers, electronic publications, and digital content.

\subsection*{Energy}
- \textbf{CVX – Chevron Corporation}: Engages in integrated energy and
chemicals operations, including exploration, production, refining, and
marketing of oil and gas products.\\   
- \textbf{XOM – Exxon Mobil Corporation}: Explores and produces crude oil and natural gas, and manufactures petroleum products, operating globally.  \\
- \textbf{BP – BP p.l.c.}: An integrated energy company providing carbon products and services, operating through Gas \& Low Carbon Energy, oil production, and customers \& products segments.  \\
- \textbf{Shanxi Coking Coal – Shanxi Coking Coal Energy Group Co., Ltd.}: Produces and sells various coal products, including coking coal, fat coal, gas coal, and lean coal in China.  \\
- \textbf{SINOPEC – China Petroleum \& Chemical Corporation}: Engages in integrated energy operations, including exploration, production, refining, and marketing of petroleum and petrochemical products.  \\
- \textbf{SHENERGY – Shenergy Company Limited}: Develops, constructs, and manages electric power, oil, and natural gas projects in China.\\

\subsection*{Technology}
- \textbf{AAPL – Apple Inc.}: Designs, manufactures, and markets smartphones, personal computers, tablets, wearables, and accessories worldwide.  \\
- \textbf{MSFT – Microsoft Corporation}: Develops and supports software, services, devices, and solutions, including operating systems, productivity applications, and cloud services.  \\
- \textbf{AMD – Advanced Micro Devices, Inc.}: Operates as a semiconductor company, offering CPUs, GPUs, and other computing solutions across data center, client, gaming, and embedded markets.  \\
- \textbf{ZTE – ZTE Corporation}: Provides integrated information and communication technology solutions, including wireless, wireline, devices, and telecommunication software systems.  \\
- \textbf{IFLYTEK – iFLYTEK CO., LTD.}: Engages in artificial intelligence technologies, offering products like smart translators, recorders, and voice-based educational tools.  \\
- \textbf{Inspur Electronic Information – Inspur Electronic Information Industry Co., Ltd.}: Provides information technology infrastructure products, including servers, storage solutions, and cloud computing services.

\subsection*{Healthcare}
- \textbf{CAH – Cardinal Health, Inc.}: Operates as a healthcare services and products company, distributing pharmaceuticals and medical products globally.  \\
- \textbf{PFE – Pfizer Inc.}: Discovers, develops, manufactures, and markets biopharmaceutical products, including vaccines and medicines across various therapeutic areas.  \\
- \textbf{BSX – Boston Scientific Corporation}: Develops, manufactures, and markets medical devices used in various interventional medical specialties worldwide.  \\
- \textbf{Yunnan Baiyao – Yunnan Baiyao Group Co., Ltd.}: Produces traditional Chinese medicines and healthcare products, including herbal remedies and personal care items.  \\
- \textbf{Zhifei Biological Products – Chongqing Zhifei Biological Products Co., Ltd.}: Develops and manufactures vaccines and biological products for disease prevention.  \\
- \textbf{Aier Eye Hospital – Aier Eye Hospital Group Co., Ltd.}: Operates a network of specialized eye hospitals and clinics, providing ophthalmic medical services.

\subsection*{Financials}
- \textbf{AMP – Ameriprise Financial, Inc.}: Provides financial planning, asset management, and insurance services to individuals and institutions.  \\
- \textbf{BAC – Bank of America Corporation}: Offers banking and financial products and services for individuals, small- and middle-market businesses, institutional investors, large corporations, and governments worldwide.  \\
- \textbf{STT – State Street Corporation}: Provides financial services and products to institutional investors worldwide, including investment servicing, investment management, and investment research and trading.  \\
- \textbf{Merchants Securities – China Merchants Securities Co., Ltd.}: Engages in securities brokerage, investment banking, asset management, and other financial services in China.  \\
- \textbf{Ping An Insurance – Ping An Insurance (Group) Company of China, Ltd.}: Provides insurance, banking, and financial services, including life and health insurance, property and casualty insurance, and asset management.  \\
- \textbf{Agricultural Bank of China – Agricultural Bank of China Limited}: Offers banking products and services to individuals, enterprises, and government agencies, including deposits, loans, and wealth management.\\

\subsection*{Consumer Discretionary}
- \textbf{AMZN – Amazon.com, Inc.}: Engages in the retail sale of consumer products and subscriptions through online and physical stores, and provides cloud computing services.  \\
- \textbf{BKNG – Booking Holdings Inc.}: Provides online travel and related services, including accommodation reservations, rental cars, and airline tickets.  \\
- \textbf{HD – The Home Depot, Inc.}: Operates as a home improvement retailer, offering building materials, home improvement products, and lawn and garden products.  \\
- \textbf{BYD – BYD Company Limited}: Engages in the manufacture and sale of automobiles, rechargeable batteries, and photovoltaic products.  \\
- \textbf{Midea – Midea Group Co., Ltd.}: Manufactures and sells home appliances, HVAC systems, and robotics and automation systems.  \\
- \textbf{Hisense Home Appliances – Hisense Home Appliances Group Co., Ltd.}: Produces and sells household electrical appliances, including refrigerators, air conditioners, and washing machines.

\subsection*{Consumer Staples}
- \textbf{COST – Costco Wholesale Corporation}: Operates membership warehouses that offer branded and private-label products across a wide range of merchandise categories.  \\
- \textbf{WMT – Walmart Inc.}: Operates retail stores, including supermarkets, discount stores, and warehouse clubs, offering a wide assortment of merchandise and services.  \\
- \textbf{TGT – Target Corporation}: Operates as a general merchandise retailer, offering food assortments, apparel, and home furnishings.  \\
- \textbf{Moutai – Kweichow Moutai Co., Ltd.}: Produces and sells Moutai liquor and other alcoholic beverages in China.  \\
- \textbf{Wuliangye – Wuliangye Yibin Co., Ltd.}: Engages in the production and sale of liquor products, including the Wuliangye series of liquors.  \\
- \textbf{Haitian Flavouring and Food – Haitian Flavouring and Food Co., Ltd.}: Produces and sells flavoring products, including soy sauce, oyster sauce, and vinegar.

\subsection*{Industrials}
- \textbf{DOV – Dover Corporation}: Provides equipment, components, consumable supplies, aftermarket parts, software and digital solutions, and support services worldwide.  \\
- \textbf{LMT – Lockheed Martin Corporation}: An aerospace and defense company, engages in the research, design, development, manufacture, integration, and sustainment of advanced technology systems, products, and services.  \\
- \textbf{EFX – Equifax Inc.}: Operates as a data, analytics, and technology company, providing information solutions and human resources business process outsourcing services.  \\
- \textbf{Shenzhen Airport – Shenzhen Airport Co., Ltd.}: Operates and manages Shenzhen Bao'an International Airport in China.  \\
- \textbf{CRRC – CRRC Corporation Limited}: Engages in the research and development, design, manufacture, refurbishment, sale, leasing, and technical services of railway transportation equipment.  \\
- \textbf{SANY HEAVY INDUS – Sany Heavy Industry Co., Ltd.}: Engages in the research and development, manufacture, and sale of construction machinery.

\subsection*{Materials}
- \textbf{LIN – Linde plc}: Operates as an industrial gas company,
\sid{providing oxygen, nitrogen, argon CO2 and hydrogen as well as
engineering and process technologies.}\\
- \textbf{MOS – The Mosaic Company}: Produces and markets concentrated phosphate and potash crop nutrients for the global agriculture industry.  \\
- \textbf{ECL – Ecolab Inc.}: Provides water, hygiene, and infection prevention solutions and services, offering comprehensive programs and services to promote safe food, maintain clean environments, and optimize water and energy use.  \\
- \textbf{BAO IRON – Baoshan Iron \& Steel Co., Ltd.}: Engages in the manufacture and sale of iron and steel products.  \\
- \textbf{Ganfeng Lithium – Ganfeng Lithium Group Co., Ltd.}: Manufactures and sells lithium products, including lithium compounds and metal, and provides lithium battery recycling services.  \\
- \textbf{WANHUA CHEM – Wanhua Chemical Group Co., Ltd.}: Provides polyurethane, petrochemical, and performance chemicals and materials.

\subsection*{Real Estate}
- \textbf{AMT – American Tower Corporation}: Owns, operates, and develops multitenant communications real estate, including wireless and broadcast towers.  \\
- \textbf{CSGP – CoStar Group, Inc.}: Provides information, analytics, and online marketplace services to the commercial real estate industry.  \\
- \textbf{DLR – Digital Realty Trust, Inc.}: Owns, acquires, develops, and operates data centers, providing colocation and interconnection solutions.  \\
- \textbf{POLY DEVELOPMENTS – Poly Developments and Holdings Group Co., Ltd.}: Engages in real estate development and property management services.  \\
- \textbf{URBAN CONS DEV – Beijing Urban Construction Investment \& Development Co., Ltd.}: Involved in urban infrastructure construction and real estate development.  \\
- \textbf{Grandjoy – Grandjoy Holdings Group Co., Ltd.}: Engages in real estate development, commercial property operation, and property management services.

\subsection*{Utilities}
- \textbf{AWK – American Water Works Company, Inc.}: Provides water and wastewater services to residential, commercial, industrial, and other customers.  \\
- \textbf{AEE – Ameren Corporation}: Generates and distributes electricity and distributes natural gas to customers in the United States.  \\
- \textbf{SRE – Sempra Energy}: Operates as an energy infrastructure company, focusing on electric and gas infrastructure and utilities.  \\
- \textbf{WUHAN HOLDING – Wuhan Holding Co., Ltd.}: Engages in water supply, sewage treatment, and other utility services.  \\
- \textbf{Dazhong Public – Dazhong Public Utilities Group Co., Ltd.}: Provides public utility services, including gas supply and sewage treatment.  \\
- \textbf{SEP – Shanghai Electric Power Co., Ltd.}: Engages in the generation and distribution of electric power.

\end{document}